\documentclass[reqno,a4paper,draft]{amsart}
\usepackage{enumitem}

\setenumerate{label=\textnormal{(\arabic*)}}

\usepackage{amsmath,amssymb,dsfont,verbatim,mathtools,bm,geometry,fge}

\usepackage[latin1]{inputenc}
\usepackage[raggedright]{titlesec}
\usepackage{mathtools}
\usepackage{tikz}
\usepackage{float,subfig}
\usepackage{amsrefs}

\titleformat{\chapter}[display]
{\normalfont\huge\bfseries}{\chaptertitlename\\thechapter}{20pt}{\Huge}
\titleformat{\section}
{\normalfont\Large\bfseries\center}{\thesection}{1em}{}
\titleformat{\subsection}
{\normalfont\large\bfseries}{\thesubsection}{1em}{}
\titleformat{\subsubsection}[runin]
{\normalfont\normalsize\bfseries}{\thesubsubsection}{1em}{}
\titleformat{\paragraph}[runin]
{\normalfont\normalsize\bfseries}{\theparagraph}{1em}{}
\titleformat{\subparagraph}[runin]
{\normalfont\normalsize\bfseries}{\thesubparagraph}{1em}{}
\titlespacing*{\chapter} {0pt}{50pt}{40pt}
\titlespacing*{\section} {0pt}{3.5ex plus 1ex minus .2ex}{2.3ex plus .2ex}
\titlespacing*{\subsection} {0pt}{3.25ex plus 1ex minus .2ex}{1.5ex plus .2ex}
\titlespacing*{\subsubsection}{0pt}{3.25ex plus 1ex minus .2ex}{1.5ex plus .2ex}
\titlespacing*{\paragraph} {0pt}{3.25ex plus 1ex minus .2ex}{1em}
\titlespacing*{\subparagraph} {\parindent}{3.25ex plus 1ex minus .2ex}{1em}
\newtheorem{theorem}{Theorem}[section]
\newtheorem{lemma}[theorem]{Lemma}

\theoremstyle{definition}

\theoremstyle{remark}

\DeclareMathOperator{\Aut}{Aut}

\DeclareMathOperator{\Jac}{Jac}

\DeclareMathOperator{\en}{en}

\DeclareMathOperator{\st}{st}

\DeclareMathOperator{\Succ}{Succ}
\DeclareMathOperator{\Pred}{Pred}
\DeclareMathOperator{\Dir}{Dir}

\newcommand{\ov}{\overline}

\renewcommand{\theequation}{\thesection.\arabic{equation}}

\begin{document}

\title{A differential equation for polynomials related to the Jacobian conjecture}

\author{Christian Valqui}
\address{Pontificia Universidad Cat\'olica del Per\'u, Secci\'on Matem\'aticas, PUCP, Av. Universitaria 1801, San Miguel, Lima 32, Per\'u.}

\address{Instituto de Matem\'atica y Ciencias Afines (IMCA) Calle Los Bi\'ologos 245. Urb San C\'esar. La Molina, Lima 12, Per\'u.}
\email{cvalqui@pucp.edu.pe}

\thanks{Christian Valqui was supported by PUCP-DGI-2013-3036, PUCP-DGI-2012-0011}

\author{Jorge A. Guccione}
\address{Departamento de Matem\'atica\\ Facultad de Ciencias Exactas y Naturales-UBA, Pabell\'on~1-Ciudad Universitaria\\ Intendente Guiraldes 2160 (C1428EGA) Buenos Aires, Argentina.}
\address{Instituto de Investigaciones Matem\'aticas ``Luis A. Santal\'o"\\ Facultad de Ciencias Exactas y Natu\-rales-UBA, Pabell\'on~1-Ciudad Universitaria\\ Intendente Guiraldes 2160 (C1428EGA) Buenos Aires, Argentina.}
\email{vander@dm.uba.ar}

\author{Juan J. Guccione}
\address{Departamento de Matem\'atica\\ Facultad de Ciencias Exactas y Naturales-UBA\\ Pabell\'on~1-Ciudad Universitaria\\ Intendente Guiraldes 2160 (C1428EGA) Buenos Aires, Argentina.}
\address{Instituto Argentino de Matem\'atica-CONICET\\ Savedra 15 3er piso\\ (C1083ACA) Buenos Aires, Ar\-gentina.}
\email{jjgucci@dm.uba.ar}

\thanks{Jorge A. Guccione and Juan J. Guccione were supported by PIP 112-200801-00900 (CONICET)}

\begin{abstract}
We analyze a possible minimal counterexample to the Jacobian Conjecture $P,Q$ with
$\gcd(\deg(P),\deg(Q))=16$ and
show that its existence depends only on the existence of solutions for a certain
Abel differential equation of the second kind.
\end{abstract}

\maketitle

{\textrm RESUMEN:} Analizamos un posible contraejemplo $P,Q$ a la conjetura del jacobiano con
$\gcd(\deg(P),\deg(Q))=16$ y mostramos que su existencia depende exclusivamente de la
existencia de soluciones de una cierta ecuación diferencial de Abel de segundo tipo.

Palabras claves: Jacobiano, ecuación diferencial de Abel.

\section{Introduction}
\renewcommand{\thetheorem}{\arabic{theorem}}
\renewcommand{\theequation}{\arabic{equation}}
In a recent article~\cite{G-G-V3}, we managed to
describe the shape of possible minimal counterexample to JC (the Jacobian conjecture as stated in~\cite{K})
given by a pair of polynomials $(P,Q)$ with $\gcd(\deg(P),\deg(Q))=B$, where
$$
B := \begin{cases}\infty & \text{if JC is true,}\\ \min\bigl(\gcd(\deg(P),\deg(Q)\bigr)&\text{where
$(P,Q)$ is a counterexample to JC, if JC is false.}
\end{cases}
$$
We arrived at the following theorem:
\begin{theorem}[\cite{G-G-V3}*{Theorem 8.10}] \label{teorema principal} If $B=16$, then there exist
$\mu_0,\mu_1,\mu_2,\mu_3\in K$ with $\mu_0\ne 0$ and $P, Q\in L:= K[x,y]$
such that
$$
P,Q\in L,\quad \ell_{1,-1}(P)=x^3 y+\mu_3 x^{2},\quad\ell_{1,-1}(Q)=x^2y+\mu_3 x
$$
and
\begin{equation}\label{Jacobiano}
\quad[P,Q]=x^4 y+\mu_0 +\mu_1 x+\mu_2 x^2+\mu_3 x^3.
\end{equation}
Moreover, there exists $j\in\mathds{N}$ such that $(j,1)\in\Dir(P)\cap\Dir(Q)$,
$$
\st_{j,1}(P)=(3,1),\quad \st_{j,1}(Q)=(2,1),\quad \en_{j,1}(P)=(0,m)\quad\text{and}\quad
\en_{j,1}(Q)=(0,n),
$$
where $m=3j+1$ and $n=2j+1$.
\end{theorem}
\noindent By \cite{H}*{Theorem 2.23} we know that $B\ge 16$. Hence, if we can prove that such a pair cannot exist, necessarily $B>16$.

\noindent In Section~1 we will
show how the existence of such a pair $(P,Q)$ would allow the construction of
a counterexample to the Jacobian Conjecture. We use the notations of~\cite{G-G-V3}.

\noindent In Section 2 we write, according to Theorem~\ref{teorema principal},
$$
P=x^3 y+x^2 p_2(y)+x p_1(y)+p_0(y)\quad\text{and}\quad Q=x^2 y+x q_1(y)+q_0(y).
$$
Then the condition~\eqref{Jacobiano} translates into a system of four first order
differential equations for the polynomials $p_0,p_1,q_0,q_1,q_2$. We reduce
this system to a single equation for two polynomials and we prove the following theorem:

\begin{theorem}
\label{teorema principal1} $B=16$ if and only if there exist $A,q_1\in  K[y]$ and
$\mu_0,\mu_1,\mu_2,\mu_3\in K$ with $\mu_0\ne 0$,
\begin{equation}\label{condiciones1}
  A(0)=-\frac 14 \mu_3^2,\quad A'(0)=\mu_2\quad\text{and}\quad \mu_3 A''(0)=-6\mu_1-2\mu_3 q_1''(0),
\end{equation}
such that
\begin{equation}\label{EDPol1}
    6\left(A-\frac{q_1^2}4+\frac{\mu_3}4 q_1-\frac{\mu_2}6 y\right)^2=
    4yAA'+6\left(\frac{\mu_3}4 q_1-\frac{\mu_2}6 y\right)^2-\mu_2 y q_1^2+3 \mu_1 y^2 q_1-6 \mu_0 y^3.
\end{equation}
\end{theorem}

We were not able to obtain a solution of~\eqref{EDPol1} satisfying~\eqref{condiciones1} with $\mu_0\ne 0$
(which would yield a counterexample to the JC), nor could we discard the existence
of such a solution (which would prove $B>16$).
We analyze some particular cases of \eqref{EDPol1}, for example we show that
for $\mu_3=\mu_2=\mu_1=\mu_0=0$ the only possible solutions are $(\rho,\sigma)$-homogeneous
for $(\rho,\sigma)=(j,1)$, where $j+1=\deg(q_1)$. We also recognize~\eqref{EDPol1} as an
Abel differential equation of second kind, for which no general solution is known. Using
a standard trick we write this equation in a shorter form in~\eqref{Ade} and in~\eqref{Ade1}.

\section{Construction of an counterexample}
\renewcommand{\thetheorem}{\thesection.\arabic{theorem}}
\renewcommand{\theequation}{\thesection.\arabic{equation}}
We reverse the order of the construction leading to Theorem 8.10 of \cite{G-G-V3}.
Starting from a pair $(P,Q)$ as in Theorem~\ref{teorema principal},
we apply different automorphisms of $L$ and $L^{(1)}$ and obtain a counterexample $(\tilde P,\tilde Q)$ with
$\gcd(\deg(\tilde P),\deg(\tilde Q))=16$.

Recall from
\cite{G-G-V3} the automorphisms $\psi_1\in \Aut(L)$ and $\psi_3\in \Aut(L^{(1)})$ given by
\begin{xalignat*}{3}
&\psi_1(x):=y,   &&\psi_3(x):=x^{-1},\\
&\psi_1(y):=-x,  &&\psi_3(y):=x^3 y.
\end{xalignat*}
For $(\rho,\sigma)\in \ov{\mathfrak{V}}$ and $k\in\{1,3\}$, we define
$(\rho_k,\sigma_k):=\ov \psi_k(\rho,\sigma)$ by
\begin{align*}
&\ov \psi_1(\rho,\sigma):=(\sigma,\rho)\quad\text{and}\quad
\ov \psi_3(\rho,\sigma):=\begin{cases} (-\rho,3\rho+\sigma)& \text{if $(\rho,\sigma)\le (-1,2)$,}\\
(\rho,-3\rho-\sigma)& \text{if $(\rho,\sigma)> (-1,2)$.}
\end{cases}
\end{align*}
We have following lemma:
\begin{lemma}[\cite{G-G-V3}*{Lemma 6.6}]\label{rho por automorphismos} Let $P\in L^{(1)}$. The maps $\psi_1$ and $\psi_3$
satisfy the following properties:
\begin{enumerate}

\smallskip

\item For all $i,j\in \mathds{N}_0$ we have
$v_{\rho_1,\sigma_1}(\psi_1(x^i y^j))=v_{\rho,\sigma}(x^i y^j)$,
 and if $P\in L$, then
$$
\ell_{\rho_1,\sigma_1}(\psi_1(P)) = \psi_1\left(\ell_{\rho,\sigma}(P)\right)\quad\text{and}\quad
\ell\ell_{\rho_1,\sigma_1}(\psi_1(P))= \psi_1\left(\ell\ell_{\rho,\sigma}(P)\right).
$$

\smallskip

\item If $(\rho,\sigma)\le (-1,2)$, then $v_{\rho_3,\sigma_3}(\psi_3(x^i y^j))=
v_{\rho,\sigma}(x^i y^j)\qquad\text{for all $i\in \mathds{N}_0$ and
$j\in \mathds{Z}$}$,
$$
\ell_{\rho_3,\sigma_3}(\psi_3(P)) = \psi_3\left(\ell_{\rho,\sigma}(P)\right)\quad\text{and}\quad
\ell\ell_{\rho_3,\sigma_3}(\psi_3(P))= \psi_3\left(\ell\ell_{\rho,\sigma}(P)\right).
$$

\smallskip

\item If $(\rho,\sigma)>(-1,2)$, then $v_{\rho_3,\sigma_3}(\psi_3(x^i y^j))=
-v_{\rho,\sigma}(x^i y^j)\qquad\text{for all $i\in \mathds{N}_0$ and
$j\in \mathds{Z}$}$,
$$
\ell_{\rho_3,\sigma_3}(\psi_3(P))=\psi_3\left(\ell\ell_{\rho,\sigma}(P)\right)\quad\text{and}\quad
\ell\ell_{\rho_3,\sigma_3}(\psi_3(P))= \psi_3\left(\ell_{\rho,\sigma}(P)\right).
$$
\end{enumerate}
\end{lemma}

Moreover clearly $\Jac(\psi_1)=[\psi_1(x),\psi_1(y)]=1$ and $\Jac(\psi_3)=-x$.\\
Let $(P,Q)$ be as in Theorem~\ref{teorema principal}.

\bigskip

\noindent{\bf FIRST STEP:}\\
Set $P_1:=\psi_3(P)$ and $Q_1:=\psi_3(Q)$ and $(\tilde\rho,\tilde\sigma):=(-j,3j+1)$.
Using Lemma~\ref{rho por automorphismos}
 one checks that
$\Pred_{P_1}(\tilde\rho,\tilde\sigma)=\Pred_{Q_1}(\tilde\rho,\tilde\sigma)=(1,-1)$,
$$
\en_{\tilde\rho,\tilde\sigma}(P_1)=(0,1),\quad \en_{\tilde\rho,\tilde\sigma}(Q_1)=(1,1),\quad
w(\ell\ell_{-1,3}( P_1))=m(3,1),\quad w(\ell\ell_{-1,3}( Q_1))=n(3,1)
$$
and
$$
\ell_{-1,2}(P_1)=y+\mu_3 x^{-2}\quad\text{and}\quad\ell_{-1,2}(Q_1)=xy+\mu_3 x^{-1},
$$
where $m:=3j+1$ and $n:=2j+1$.
Moreover, using that
$$
[\varphi(P),\varphi(Q)]=\varphi([P,Q])[\varphi(x),\varphi(y)],
$$
for all morphisms $\varphi$, we obtain
$[P_1,Q_1]=-(y+\mu_0x+\mu_1+\mu_2 x^{-1}+\mu_3 x^{-2})$.
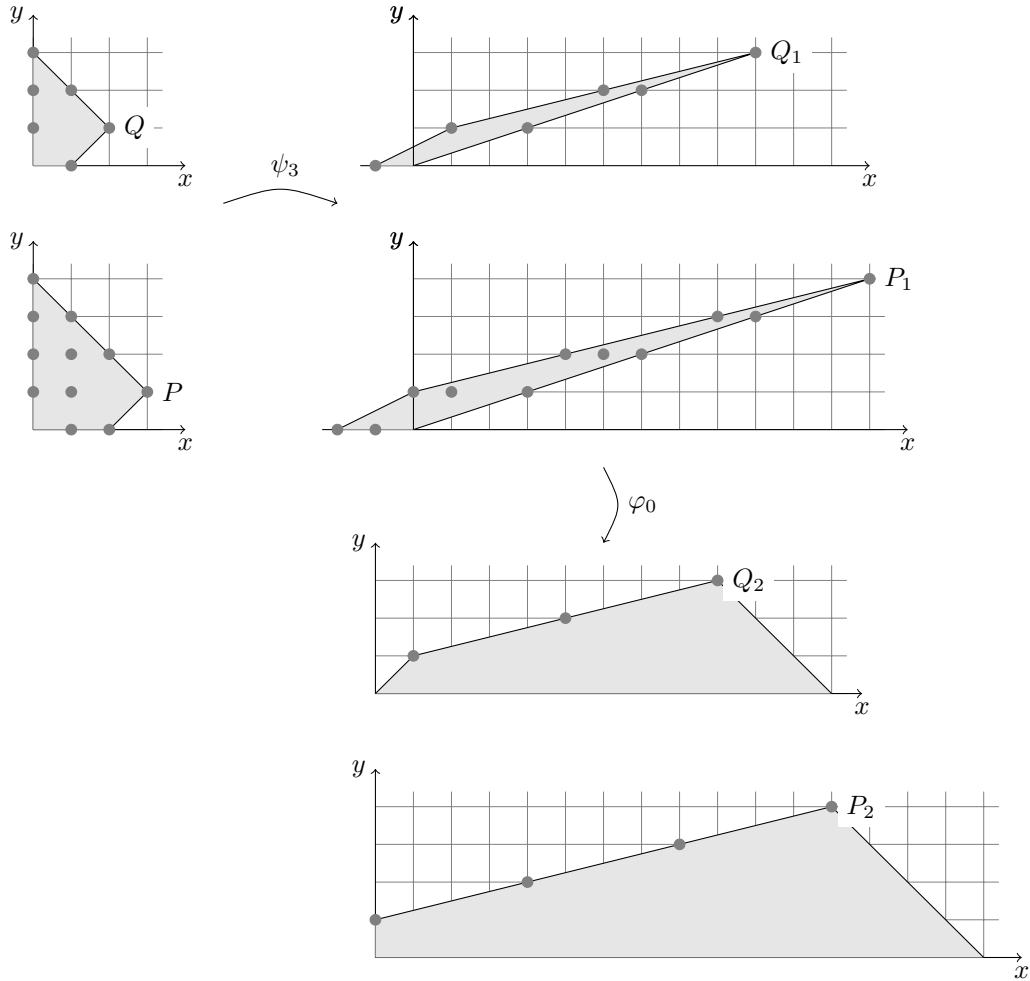
\begin{figure}[ht]
\centering
\begin{tikzpicture}
\draw[step=.5cm,gray,very thin] (0,3.5) grid (1.7,5.2);
\draw [->] (0,3.5) -- (2,3.5) node[anchor=north]{$x$};
\draw [->] (0,3.5) --  (0,5.5) node[anchor=east]{$y$};
\draw[step=.5cm,gray,very thin] (0,0) grid (1.7,2.2);
\draw [->] (0,0) -- (2,0) node[anchor=north]{$x$};
\draw [->] (0,0) --  (0,2.5) node[anchor=east]{$y$};
\draw[step=.5cm,gray,very thin] (5,3.5) grid (10.7,5.2);
\draw [->] (4.3,3.5) -- (11,3.5) node[anchor=north]{$x$};
\draw [->] (5,3.5) --  (5,5.5) node[anchor=east]{$y$};
\draw[step=.5cm,gray,very thin] (5,0) grid (11.2,2.2);
\draw [->] (3.8,0) -- (11.5,0) node[anchor=north]{$x$};
\draw [->] (5,0) --  (5,2.5) node[anchor=east]{$y$};
\fill[gray!20] (0,3.5) -- (0.5,3.5) -- (1,4)--(0,5);
\fill[gray!20] (0,0) -- (1,0) -- (1.5,0.5) -- (0,2);
\fill[gray!20] (5,3.5) -- (9.5,5) -- (5.5,4)--(4.5,3.5);
\fill[gray!20] (5,0) -- (11,2) -- (5,0.5)--(4,0);
\draw [->] (5,3.5) --  (5,5.5) node[anchor=east]{$y$};
\draw [->] (5,0) --  (5,2.5) node[anchor=east]{$y$};
\draw (1,0) --  (1.5,0.5) node[fill=white,right=2pt]{$P$} -- (0,2);
\draw (0.5,3.5) --  (1,4) node[fill=white,right=2pt]{$Q$} -- (0,5);
\draw (5,0) --  (11,2) node[fill=white,right=2pt]{$P_1$} -- (5,0.5)--(4,0);
\draw (5,3.5) --  (9.5,5) node[fill=white,right=2pt]{$Q_1$} -- (5.5,4)--(4.5,3.5);
\filldraw [gray]  (0.5,0)    circle (2pt)
                  (1,0)      circle (2pt)
                  (0,0.5)    circle (2pt)
                  (0.5,0.5)  circle (2pt)
                  (1.5,0.5)  circle (2pt)
                  (0,1)    circle (2pt)
                  (0.5,1)      circle (2pt)
                  (1,1)      circle (2pt)
                  (0,1.5)    circle (2pt)
                  (0.5,1.5)  circle (2pt)
                  (0,2)    circle (2pt)
                  (0.5,3.5)    circle (2pt)
                  (0,4)    circle (2pt)
                  (1,4)      circle (2pt)
                  (0,4.5)    circle (2pt)
                  (0.5,4.5)    circle (2pt)
                  (0,5)      circle (2pt);
\filldraw [gray]  (5.5,0.5)    circle (2pt)
                  (4.5,0)      circle (2pt)
                  (4,0)    circle (2pt)
                  (5,0.5)  circle (2pt)
                  (6.5,0.5)  circle (2pt)
                  (7,1)    circle (2pt)
                  (8,1)      circle (2pt)
                  (7.5,1)      circle (2pt)
                  (9,1.5)    circle (2pt)
                  (9.5,1.5)  circle (2pt)
                  (11,2)    circle (2pt)
                  (4.5,3.5)    circle (2pt)
                  (5.5,4)    circle (2pt)
                  (6.5,4)      circle (2pt)
                  (7.5,4.5)    circle (2pt)
                  (8,4.5)    circle (2pt)
                  (9.5,5)      circle (2pt);
\draw[->] (2.5,3) .. controls (3.25,3.25) .. (4,3);
\draw (3,3.5) node[right,text width=7cm]{$\psi_3$};
%
%lo de abajo
\draw (7.7,-1) node[right,text width=1cm]{$\varphi_0$};
\draw[->] (7.5,-0.5) .. controls (7.75,-1) .. (7.5,-1.5);
\draw[step=.5cm,gray,very thin] (4.5,-3.5) grid (10.7,-1.8);
\draw [->] (4.5,-3.5) -- (10.9,-3.5) node[anchor=north]{$x$};
\draw [->] (4.5,-3.5) --  (4.5,-1.5) node[anchor=east]{$y$};
\fill[gray!20] (4.5,-3.5) -- (10.5,-3.5) -- (9,-2)--(5,-3);
\draw (10.5,-3.5) --  (9,-2) node[fill=white,right=2pt]{$Q_2$} -- (5,-3)--(4.5,-3.5);
\draw[step=.5cm,gray,very thin] (4.5,-7) grid (12.7,-4.8);
\draw [->] (4.5,-7) -- (13,-7) node[anchor=north]{$x$};
\draw [->] (4.5,-7) --  (4.5,-4.5) node[anchor=east]{$y$};
\fill[gray!20] (4.5,-7) -- (12.5,-7) -- (10.5,-5)--(4.5,-6.5);
\draw (12.5,-7) --  (10.5,-5) node[fill=white,right=2pt]{$P_2$} -- (4.5,-6.5);
\filldraw [gray]  (9,-2)    circle (2pt)
                  (7,-2.5)      circle (2pt)
                  (5,-3)    circle (2pt)
                  (10.5,-5)  circle (2pt)
                  (8.5,-5.5)  circle (2pt)
                  (6.5,-6)    circle (2pt)
                  (4.5,-6.5)     circle (2pt);
\end{tikzpicture}
\caption{Illustration of the first two steps, for $j=1$.}
\end{figure}

\noindent{\bf SECOND STEP}\\
Set $P_2:=\varphi_0(P_1)$ and $Q_2:=\varphi_0(Q_1)$, where
$\varphi_0(y):=y-(\mu_0x+\mu_1+\mu_2 x^{-1}+\mu_3 x^{-2})$ and
$\varphi_0(x):=x$ (note that $\Jac(\varphi_0)=1$).
Then $P_2,Q_2\in L$ and
$$
[P_2,Q_2]=-y,\quad \Dir(P_2)=\Dir(Q_2)=\{(\tilde\rho,\tilde\sigma),(1,1)\},\quad
\en_{\tilde\rho,\tilde\sigma}(P_2)=(0,1),\quad \en_{\tilde\rho,\tilde\sigma}(Q_2)=(1,1)
$$
and
$$
\ell_{1,1}(P_2)=\lambda_P R_2^m\quad\text{and}\quad\ell_{1,1}(Q_2)=\lambda_Q R_2^n,
$$
for $R_2=x^3(y-\mu_0 x)$.

\bigskip

\noindent{\bf THIRD STEP}\\
Since $P_2,Q_2\in L$, we can apply $\psi_1$. We set
$P_3:=\psi_1(P_2)$, $Q_3:=\psi_1(Q_2)$ and
$(\ov\rho,\ov\sigma):=(3j+1,-j)$. Then
$$
[P_3,Q_3]=-x,\quad \Dir(P_3)=\Dir(Q_3)=\{(\ov\rho,\ov\sigma),(1,1)\},\quad
\en_{\ov\rho,\ov\sigma}(P_3)=(1,0),\quad \en_{\ov\rho,\ov\sigma}(Q_3)=(1,1)
$$
and
$$
\ell_{1,1}(P_3)=\tilde \lambda_P R_3^m\quad\text{and}\quad\ell_{1,1}(Q_3)=\tilde \lambda_Q R_3^n,
$$
for $R_3=y^3(y+\frac{1}{\mu_0}x)$.
\bigskip
\begin{figure}[h]
\centering
\begin{tikzpicture}
\draw[step=.125cm,gray,very thin] (0,3.5) grid (1.7,5.2);
\draw [->] (0,3.5) -- (2,3.5) node[anchor=north]{$x$};
\draw [->] (0,3.5) --  (0,5.5) node[anchor=east]{$y$};
\draw[step=.125cm,gray,very thin] (0,0) grid (1.7,2.2);
\draw [->] (0,0) -- (2,0) node[anchor=north]{$x$};
\draw [->] (0,0) --  (0,2.5) node[anchor=east]{$y$};
\draw[step=.125cm,gray,very thin] (5,3.5) grid (10.7,5.2);
\draw [->] (4.3,3.5) -- (11,3.5) node[anchor=north]{$x$};
\draw [->] (5,3.5) --  (5,5.5) node[anchor=east]{$y$};
\draw[step=.125cm,gray,very thin] (5,0) grid (11.2,2.2);
\draw [->] (3.8,0) -- (11.5,0) node[anchor=north]{$x$};
\draw [->] (5,0) --  (5,2.5) node[anchor=east]{$y$};
\fill[gray!20] (0,3.5) -- (0.125,3.625) -- (0.375,4.625)--(0,5);
\fill[gray!20] (0,0) -- (0.125,0) -- (0.5,1.5) -- (0,2);
\fill[gray!20] (5,3.5) -- (9.5,5) -- (8,4.625)--(5.25,3.625);
\fill[gray!20] (5,0) -- (11,2) -- (9,1.5)--(4.875,0);
\draw (0.125,0) --  (0.5,1.5) node[fill=white,right=2pt]{$P_3$} -- (0,2);
\draw (0,3.5) -- (0.125,3.625) -- (0.375,4.625) node[fill=white,right=2pt]{$Q_3$} -- (0,5);
\draw (5,0) --  (11,2) node[fill=white,right=2pt]{$P_4$} -- (9,1.5)--(4.875,0);
\draw (5,3.5) --  (9.5,5) node[fill=white,right=2pt]{$Q_4$} -- (8,4.625)--(5.25,3.625)--(5,3.5);
\filldraw [gray]  (0.125,0)    circle (1pt)
                  (0.25,0.5)      circle (1pt)
                  (0.375,1)      circle (1pt)
                  (0.5,1.5)      circle (1pt)
                  (0,2)    circle (1pt)
                  (0,3.5)  circle (1pt)
                  (0.125,3.625)  circle (1pt)
                  (0.25,4.125)      circle (1pt)
                  (0.375,4.625)    circle (1pt)
                  (0,5)      circle (1pt)
                  (5,0)      circle (1pt)
                  (11,2)    circle (1pt)
                  (9,1.5)  circle (1pt)
                  (7.625,1)  circle (1pt)
                  (6.25,0.5)  circle (1pt)
                  (4.875,0)    circle (1pt)
                  (5,3.5)    circle (1pt)
                  (9.5,5)    circle (1pt)
                  (8,4.625)      circle (1pt)
                  (6.625,4.125)  circle (1pt)
                  (5.25,3.625)    circle (1pt);
\draw[->] (2.5,3) .. controls (3.25,3.25) .. (4,3);
\draw (3,3.5) node[right,text width=7cm]{$\psi_3$};
\end{tikzpicture}
\caption{Illustration of the fourth step.}
\end{figure}
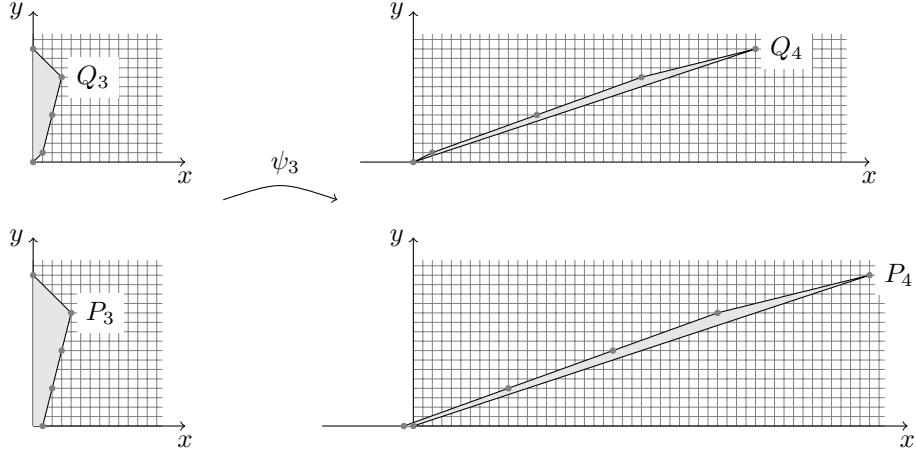

\noindent{\bf FOURTH STEP}(Figure 2)\\
We set
$P_4:=\psi_3(P_3)$, $Q_4:=\psi_3(Q_3)$ and
$(\hat\rho,\hat\sigma):=(-3j-1,8j+3)$.
Then
$$
[P_4,Q_4]=1,\quad \Dir(P_4)=\Dir(Q_4)=\{(\hat\rho,\hat\sigma),(-1,4)\},\quad
\en_{\hat\rho,\hat\sigma}(P_4)=(-1,0),\quad \en_{\hat\rho,\hat\sigma}(Q_4)=(2,1)
$$
and
$$
\ell_{-1,4}(P_4)=\tilde \lambda_P R_4^m\quad\text{and}\quad\ell_{-1,4}(Q_4)=\tilde \lambda_Q R_4^n,
$$
for $R_4=y^3 x^{12}(y+\frac{1}{\mu_0}x^{-4})$.

\bigskip

\noindent{\bf FIFTH STEP}\\
Set $P_5:=\varphi_1(P_4)$ and $Q_5:=\varphi_1(Q_4)$, where
$\varphi_1(y):=y-\frac{1}{\mu_0}x^{-4}$ and $\varphi_1(x):=x$ (note that $\Jac(\varphi_1)=1$).
Then
$$
\ell_{-1,4}(P_5)=\tilde \lambda_P R_5^m\quad\text{and}\quad\ell_{-1,4}(Q_5)=\tilde \lambda_Q  R_5^n,
$$
for $R_5=y x^{12}(y-\frac{1}{\mu_0}x^{-4})^3$.

\bigskip

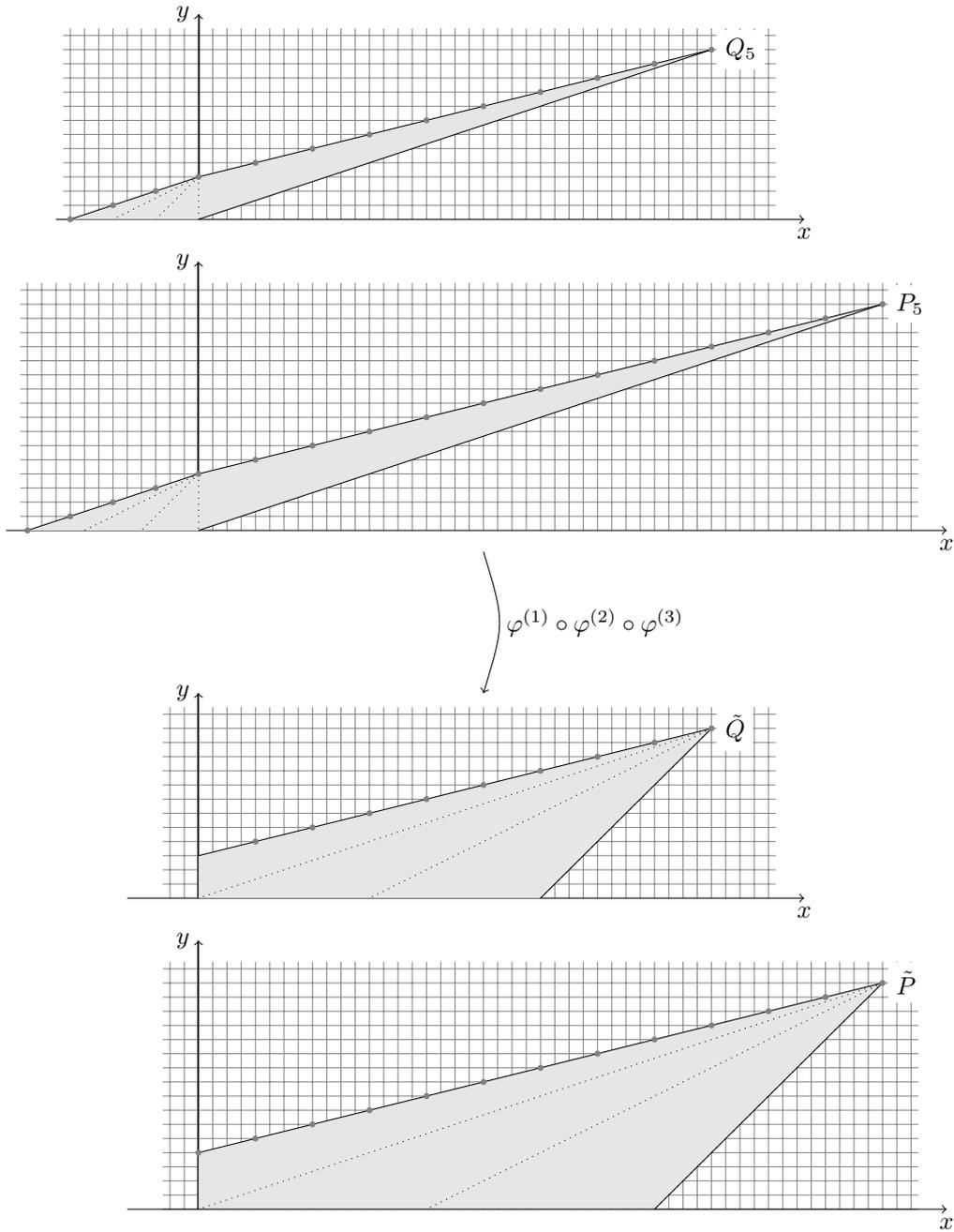
\begin{figure}[h]
\centering
\begin{tikzpicture}
\draw[step=.2cm,gray,very thin] (-0.5,0) grid (10.1,3.5);
\draw [->] (-1,0) -- (10.5,0) node[anchor=north]{$x$};
\draw [->] (0,0) --  (0,3.8) node[anchor=east]{$y$};
\draw [-] (-0.01,0) --  (-0.01,3.8) node[anchor=east]{};
\fill[gray!20] (0,0.01) -- (6.4,0.01) -- (9.6,3.2)--(0,0.8);
\draw (6.4,0) --  (9.6,3.2) node[fill=white,right=2pt]{$\tilde P$} -- (0,0.8);
\draw[dotted] (9.6,3.2) -- (0,0);
\draw[dotted] (9.6,3.2) -- (3.2,0);
\draw[dotted] (9.6,3.2) -- (6.4,0);
\draw[step=.2cm,gray,very thin] (-2.5,9.6) grid (10.1,13.1);
\draw [->] (-2.7,9.6) -- (10.5,9.6) node[anchor=north]{$x$};
\draw [->] (0,9.6) --  (0,13.4) node[anchor=east]{$y$};
\draw [-] (-0.01,9.6) --  (-0.01,13.4) node[anchor=east]{};
\fill[gray!20] (0,9.6) -- (9.6,12.8)--(0,10.4)--(-2.4,9.6);
\draw (0,9.6) --  (9.6,12.8) node[fill=white,right=2pt]{$P_5$} -- (0,10.4)--(-2.4,9.6);
\draw[dotted] (0,10.4) -- (-1.6,9.6);
\draw[dotted] (0,10.4) -- (-0.8,9.6);
\draw[dotted] (0,10.4) -- (0,9.6);
\draw[step=.2cm,gray,very thin] (-0.5,4.4) grid (8.1,7.1);
\draw [->] (-1,4.4) -- (8.5,4.4) node[anchor=north]{$x$};
\draw [->] (0,4.4) --  (0,7.3) node[anchor=east]{$y$};
\draw [-] (-0.01,4.4) --  (-0.01,7.3) node[anchor=east]{};
\fill[gray!20] (0,4.4) -- (4.8,4.4) -- (7.2,6.8)--(0,5);
\draw (4.8,4.4) --  (7.2,6.8) node[fill=white,right=2pt]{$\tilde Q$} -- (0,5);
\draw[dotted] (7.2,6.8) -- (0,4.4);
\draw[dotted] (7.2,6.8) -- (2.4,4.4);
\draw[dotted] (7.2,6.8) -- (4.8,4.4);
\draw[step=.2cm,gray,very thin] (-1.9,14) grid (8.1,16.7);
\draw [->] (-2,14) -- (8.5,14) node[anchor=north]{$x$};
\draw [->] (0,14) --  (0,16.9) node[anchor=east]{$y$};
\draw [-] (0.01,14) --  (0.01,16.9) node[anchor=east]{};
\fill[gray!20] (0,14) -- (7.2,16.4)--(0,14.6)--(-1.8,14);
\draw (0,14) --  (7.2,16.4) node[fill=white,right=2pt]{$Q_5$} -- (0,14.6)--(-1.8,14);
\draw[dotted] (0,14.6) -- (-1.2,14);
\draw[dotted] (0,14.6) -- (-0.6,14);
\draw[dotted] (0,14.6) -- (0,14);
\filldraw [gray]  (0,0.8)    circle (1pt)
                  (0.8,1)      circle (1pt)
                  (1.6,1.2)      circle (1pt)
                  (2.4,1.4)      circle (1pt)
                  (3.2,1.6)    circle (1pt)
                  (4,1.8)  circle (1pt)
                  (4.8,2)  circle (1pt)
                  (5.6,2.2)      circle (1pt)
                  (6.4,2.4)    circle (1pt)
                  (7.2,2.6)      circle (1pt)
                  (8,2.8)      circle (1pt)
                  (8.8,3)    circle (1pt)
                  (9.6,3.2)  circle (1pt)
                  (0,10.4)    circle (1pt)
                  (0.8,10.6)      circle (1pt)
                  (1.6,10.8)      circle (1pt)
                  (2.4,11)      circle (1pt)
                  (3.2,11.2)    circle (1pt)
                  (4,11.4)  circle (1pt)
                  (4.8,11.6)  circle (1pt)
                  (5.6,11.8)      circle (1pt)
                  (6.4,12)    circle (1pt)
                  (7.2,12.2)      circle (1pt)
                  (8,12.4)      circle (1pt)
                  (8.8,12.6)    circle (1pt)
                  (9.6,12.8)  circle (1pt)
                  (-0.6,10.2)    circle (1pt)
                  (-1.2,10)      circle (1pt)
                  (-1.8,9.8)      circle (1pt)
                  (-2.4,9.6)      circle (1pt);
\filldraw [gray]  (7.2,16.4)     circle (1pt)
                  (6.4,16.2)      circle (1pt)
                  (5.6,16)      circle (1pt)
                  (4.8,15.8)      circle (1pt)
                  (4,15.6)    circle (1pt)
                  (3.2,15.4)  circle (1pt)
                  (2.4,15.2)  circle (1pt)
                  (1.6,15)      circle (1pt)
                  (0.8,14.8)    circle (1pt)
                  (0,14.6)      circle (1pt)
                  (-0.6,14.4)      circle (1pt)
                  (-1.2,14.2)    circle (1pt)
                  (-1.8,14)  circle (1pt)
                  (7.2,6.8)     circle (1pt)
                  (6.4,6.6)      circle (1pt)
                  (5.6,6.4)      circle (1pt)
                  (4.8,6.2)      circle (1pt)
                  (4,6)    circle (1pt)
                  (3.2,5.8)  circle (1pt)
                  (2.4,5.6)  circle (1pt)
                  (1.6,5.4)      circle (1pt)
                  (0.8,5.2)    circle (1pt)
                  (0,5);
\draw (4.2,8.3) node[right,text width=3cm]{$\varphi^{(1)}\circ \varphi^{(2)}\circ \varphi^{(3)}$};
\draw[->] (4,9.3) .. controls (4.3,8.3) .. (4,7.3);
\end{tikzpicture}
\caption{Illustration of the sixth step.}
\end{figure}
\noindent{\bf SIXTH STEP}(Figure 3)\\
If $P_5,Q_5\in L$, then we have a counterexample to JC, since $[P_5,Q_5]=1$, $\deg(P)=16 m$ and
$\deg(Q)=16 n$ with $m\nmid n$ and $n\nmid m$.

Else set $(\rho_1,\sigma_1):=\Succ_{P_5}(-1,4)$. Then
$[\ell_{\rho_1,\sigma_1}(P_5),\ell_{\rho_1,\sigma_1}(Q_5)]=0$ and so
$$
\ell_{\rho_1,\sigma_1}(P_5)=\hat\lambda_P R_6^m\quad\text{and}\quad\ell_{-1,4}(Q_5)=\hat \lambda_Q R_6^n,
$$
for some $R_6=y+\lambda_k x^{-k}$ with $\lambda_k\in K^{\times}$ and $k\in\{1,2,3\}$. Note that
$(\rho_1,\sigma_1)=(-1,k)$. If necessary, we apply successively $\varphi^{(k)}$ given by
$\varphi^{(k)}(y):=y-\lambda_k x^{-k}$ and $\varphi^{(k)}(x):=x$, to obtain finally
the desired counterexample $(\tilde P,\tilde Q)$ given by
$$
(\tilde P,\tilde Q):=(\varphi^{(1)}(\varphi^{(2)}(\varphi^{(3)}(P_5))),
\varphi^{(1)}(\varphi^{(2)}(\varphi^{(3)}(Q_5))))\in L.
$$
\section{Differential equations for polynomials}

According to Theorem~\ref{teorema principal} we write
$$
P=x^3 y+x^2 p_2(y)+ xp_1(y)+p_0(y)\quad\text{and}\quad Q=x^2 y+x q_1(y)+ q_0(y).
$$
Then the equality~\eqref{Jacobiano} yields
\begin{align*}
    x^4 y&=[x^3 y,x^2 y]\\
\mu_3 x^3&=[x^3 y,x q_1(y)]+[x^2 p_2(y),x^2 y]\\
\mu_2 x^2&=[x^3 y, q_0(y)]+[x^2 p_2(y),x q_1(y)]+[x p_1(y), x^2 y]\\
\mu_1 x&=[x^2 p_2(y), q_0(y)]+[x p_1(y),x q_1(y)]+[p_0(y), x^2 y]\\
\mu_0&= [x p_1(y), q_0(y)]+[ p_0(y),x q_1(y)].
\end{align*}
The first equality is trivially true. Noting that
$$
[x^k p_k(y),x^j q_j(y)]=x^{k+j-1}(k p_k(y)q_j'(y)-jp_k'(y)q_j(y) ),
$$
we obtain the system of four differential equations for the five polynomials
$p_0,p_1,p_2,q_0,q_1$:
\begin{align*}
\mu_3 &=3 y q_1'-q_1+2 p_2-2yp_2'\\
\mu_2 &=3 y q_0'+2p_2 q_1'-p_2'q_1+p_1-2yp_1'\\
\mu_1 &=2p_2q_0'+p_1 q_1'-p_1'q_1-2yp_0'\\
\mu_0&= p_1 q_0'-p_0'q_1.
\end{align*}
Note that $\ell_{1,-1}(P)=x^3 y+\mu_3 x^{2}$ and $\ell_{1,-1}(Q)=x^2y+\mu_3 x$
imply $q_1(0)=\mu_3$ and $p_2(0)=\mu_3$. Moreover,
if we write $P=\sum_{i,j}a_{i,j}x^i y^j$, then we can assume $ a_{2,1}=p'_2(0)=0$,
replacing $P$ by $P-a_{2,1}Q$. Writing $Q=\sum_{i,j}b_{i,j}x^i y^j$, $[P,Q]=\sum_{i,j}c_{i,j}x^i y^j$
and noting that
\begin{equation}\label{coeficientes jacobiano}
    c_{i,j}=\sum_{(k,l)+(s,t)=(i,j)+(1,1)}(kt-ls)a_{k,l}b_{s,t},
\end{equation}
 one verifies that
$$
0=c_{3,1}=2 a_{3,1}b_{1,1}=2b_{1,1},
$$
using $b_{2,0}=b_{2,2}=a_{3,2}=a_{3,0}=0$ and $b_{3,k}=a_{4,k}=0$ for all $k$.
It follows that $q'_1(0)=b_{1,1}=0$ and so we can and will assume
$$
q_1(0)=\mu_3,\quad q'_1(0)=0, \quad p_2(0)=\mu_3\quad\text{and}\quad p'_2(0)=0.
$$
This allows to solve the first equation in full generality. In fact,
write $q_1=\mu_3+y^2 F'$ and $p_2=\mu_3+y G$ for some $F,G\in K[y]$. From the first equation
we obtain
$$
\mu_3=3y(2y F'+y^2 F'')-(\mu_3+y^2 F')+2(\mu_3+y G)-2y(G+y G'),
$$
from which we deduce the equality
$$
2G'=5F'+3yF''=(2F+3yF')'
$$
and so $G=F+(3/2)yF'+const$. Since $G(0)=0$, we can assume $F(0)=0$ and $G=F+(3/2)yF'$.
Hence the general solution to the first equation is
$q_1=\mu_3+y^2 F'$ and $p_2=\mu_3+y F+(3/2)y^2 F'$, for any choice of $F\in yK[y]$.

Using the second equation we can express $q_0'$ as a function of $F$ and $p_1$:
\begin{equation}\label{q0prima}
q_0'=\frac{-2 p_1 + 2 \mu_2 + 2  \mu_3 F + 4 y p_1'  - 6 y^2 F F'  -  \mu_3 y^2 F'' -
    4  y^3 (F')^2 - 4  y^3 F F'' - 3  y^4 F' F''}{6y}
\end{equation}
The third equation yields
$p_0'$ as a function of $F, p_1$ and $q_0'$:
\begin{equation}\label{p0prima}
p_0'= \frac{y p_1  (2 F' + y F'' )- \mu_1
     -  p_1'  (\mu_3 +  y^2 F')+(2 \mu_3 + y (2 F + 3 y F' ))q_0'}{2y}
\end{equation}
 Inserting the values into the fourth equation we obtain a (very big) differential equation for
$p_1$ and $F$:
\begin{equation}\label{ecuacion central}
\begin{array}{cc} 6 \mu_0 y^2\hspace{-1pt} &= y p_1
\Bigl(2( p_1 -  \mu_2 -  \mu_3 F) - 4 y p_1'  +  y^2 (6 F F'  +  \mu_3  F'') +
    4 y^3 ( (F')^2 +   F F'')  + 3 y^4 F' F'' \Bigr) \\
    &- (\mu_3 +  y^2 F') \Biggl(3  y^2 p_1 (2 F' + y F'' )-3 \mu_1 y
     - 3 y p_1'  (\mu_3 +  y^2 F')-\frac 12 (2 \mu_3 + y (2 F + 3 y F')) \\
     & \cdot(2 p_1 - 2 \mu_2 - 2  \mu_3 F - 4 y p_1'  + 6 y^2 F F'  +  \mu_3 y^2 F'' +
       4 y^3 (F')^2 + 4 y^3 F F''  + 3 y^4 F' F'' ) \Biggr)\\
       \end{array}
\end{equation}
Now we set
$$
A:=y p_1 - q_1 p_2 + \frac 34 q_1^2=-\frac 14 \mu_3^2 + yp_1  -  \mu_3 y F -
\mu_3 y^2 F' - y^3 F F' -  \frac 34 y^4(F')^2
$$
and we can express~\eqref{ecuacion central} as a differential equation for $A$ and $q_1$:
\begin{equation}\label{EDPol}
    6\left(A-\frac{q_1^2}4+\frac{\mu_3}4 q_1-\frac{\mu_2}6 y\right)^2=
    4yAA'+6\left(\frac{\mu_3}4 q_1-\frac{\mu_2}6 y\right)^2-\mu_2 y q_1^2+3 \mu_1 y^2 q_1-6 \mu_0 y^3
\end{equation}
Moreover we have
\begin{equation}\label{condiciones}
  A(0)=-\frac 14 \mu_3^2,\quad A'(0)=\mu_2\quad\text{and}\quad u_3 A''(0)=-6\mu_1-2\mu_3 q_1''(0).
\end{equation}
In fact, from the definition of $A$ we have that $A(0)=- q_1(0) p_2(0) + \frac 34 q_1(0)^2=-\frac 14 \mu_3^2$.
The other two conditions follow from the requirement that $q_0'(y)$ and $p_0'(y)$ defined
by~\eqref{q0prima} and~\eqref{p0prima} are polynomials.

This proves Theorem~\ref{teorema principal1} and is a great simplification
with respect to~\eqref{ecuacion central}, not only in the number of terms involved, but
in the type of differential equation. In fact,~\eqref{ecuacion central} is
a quadratic first order differential equation for
$A$, called an Abel differential equation of second kind. For $q_1$ it is a cuartic equation
with no derivative of $q_1$ involved.
However we were not able to obtain a solution of~\eqref{EDPol} with $\mu_0\ne 0$ and such
that~\eqref{condiciones} is satisfied (which would yield a counterexample to the JC),
nor could we discard the existence of such a solution (which would prove $B>16$).
In the sequel, we will analyze some aspects of this differential equation.

\subsection{Solutions without~{\bf \eqref{condiciones}}.}\label{seccion}
If we don't require~\eqref{condiciones}, then there exist solutions of~\eqref{EDPol}
with $\mu_0\ne 0$. Take for example $A=1-y^3-y^6/4$ and $q_1(y)=y^3+2$. Then~\eqref{EDPol}
is satisfied for $\mu_0=1$, $\mu_1=0=\mu_2$ and $\mu_3=2$. If we try to construct
a counterexample, we obtain $p_1(y)=y^5+2y^2+\frac{2}{y}\notin K[y]$. In fact this solution yields
$$
P=x^3 y+2x^2( y^3+1)+x\left(y^5 + 2 y^2 + \frac{2}{y}\right)+ \frac{y^7}7+ \frac{y^4}2 +\frac 1{y^2}\quad\text{and}
\quad Q=x^2 y+x(y^3+2)+\frac{y^5}5+y^2+\frac 2y.
$$
Note that $P,Q\in K[x,y,y^{-1}]$ and $[P,Q]=x^4 y+\mu_0 +\mu_1 x+\mu_2 x^2+\mu_3 x^3$, with $\mu_0=1\ne 0$.

\subsection{The case $\mu_3=\mu_2=\mu_1=\mu_0=0$: Homogeneous solutions.}
Consider the case $\mu_3=\mu_2=\mu_1=\mu_0=0$. Then~\eqref{EDPol} reads
$$
6\left(A-\frac{q_1^2}4\right)^2= 4yAA',
$$
and clearly, any irreducible factor of $A$ must be $y$, since any other
linear factor of $A$ would have multiplicity $2t$ on the left hand side
and $2k-1$ on the right hand side. Then we can assume $A=y^k$ for some $k$
and necessarily $q_1^2=4y^k\left(1\pm \sqrt{\frac{2k}3}\right)$, hence $k=2(j+1)$
and $q_1=2Ry^{j+1}$, for $R:=\pm 2\sqrt{1\pm \sqrt{\frac{4j+2}3}}$. Then it is straightforward to verify that
$p_2=\left(\frac 32+\frac 1{j}\right)R y^{j+1}$ and $p_1=y^{2j+1}\left(1-\left(\frac 1{j}+\frac 34 \right)R^2\right)$.
We also obtain $q_0=\lambda y^{2j+1}$ and $p_0=\lambda_1 y^{3j+1}$ for some $\lambda,\lambda_1$. Hence
$P$ and $Q$ are $(\rho,\sigma)$-homogeneous
for $(\rho,\sigma)=(j,1)$.

\subsection{Standard methods for solving Abel differential equations.}
For Abel differential equations no general solution is known. However, some methods are available:
The standard method for simplifying an Abel differential equation of the second kind
suggests the substitution $A= y^{3/2}T$ in~\eqref{EDPol}. This yields the equation
\begin{equation}\label{Ade}
    TT'=F_1(y) T+F_0(y)
\end{equation}
with
$$
F_1(y)=-\frac{1}{4 y^{5/2}}(3 q_1^2-3 \mu_3 q_1+2\mu_2 y )
$$
and
$$
F_0(y)=\frac{3}{32 y^4}( q_1^4 -2  \mu_3 q_1^3 + 4  \mu_2 y q_1^2- 8  \mu_1 y^2 q_1 + 16 \mu_0 y^3)
$$

We could't bring the equation~\eqref{Ade} into any of the 80 solvable cases listed in \cite{PZ}*{1.3.3},
nor could we discard the existence of solutions.

Following the book~\cite{PZ} we set $U=\frac 1T$ and then~\eqref{Ade} reads
\begin{equation}\label{Ade1}
    U'+F_1(y) U^2+F_0(y)U^3=0,
\end{equation}
an Abel differential equation of the first kind. Again, we couldn't find a solvable case in
~\cite{PZ} that corresponds to~\eqref{Ade1} and it is also impossible to choose $\mu_0\ne 0$,
$\mu_1$, $\mu_2$, $\mu_3$, $q_1$ and $\alpha$ such that
$$
\left(\frac{F_0}{F_1}\right)'=\alpha F_1,
$$
which is one of the known cases that allow further simplification of equation~\eqref{Ade1}.

\subsection{The case $\mu_3=0=\mu_1$.}
Let us analyze the equation~\eqref{Ade} in one particular case. Note that
by~\eqref{coeficientes jacobiano} we have
$$
\mu_1=c_{1,0}=2 a_{2,0}b_{0,1}-a_{1,1}b_{1,1}=\mu_3(2 b_{0,1}-a_{1,1}).
$$
Consequently, if $\mu_3=0$, then $\mu_1=0$. We will consider the case
$\mu_3=0=\mu_1$. In this case
$$
F_1(y)=-\frac{1}{4 y^{5/2}}(3 q_1^2-2\mu_2 y )
$$
and
$$
F_0(y)=\frac{1}{32 y^4}(3 q_1^4  + 4  \mu_2 y q_1^2+ 48 \mu_0 y^3).
$$
Again, we were unable to transform~\eqref{Ade} into one of the solvable cases of~\cite{PZ}.

We also can try to solve the case $\mu_1=0$ and $\mu_3=0$ directly in~\eqref{EDPol}.
In that case we can set $S:=\frac{q_1^2}4+\frac{\mu_2 y}6$ and then~\eqref{EDPol} reads
$$
3(A-S)^2=2y AA'-2\mu_2 y S+\frac{5}{12}\mu_2^2 y^2-3\mu_0 y^3.
$$
We couldn't find solutions with $\mu_0\ne 0$ such that $S-\frac{\mu_2 y}6$ is a square.

\subsection{Low degree cases.}
Finally we solve~\eqref{EDPol} with the initial conditions~\eqref{condiciones}
for some low degree cases. One can show that $\deg(A)=2\deg(q_1)$, and we were able to solve the cases
$\deg(q_1)=2,3,4$, assuming $q_1$ monic and setting $\mu_0,\mu_1,\mu_2,\mu_3$ and the coefficients of $q_1$ and $A$ as variables.
For $\deg(q_1)=3$ we obtain the solution $\mu_2=\mu_1=\mu_0=0$ and $A=-y^6/4-\mu_3 y^3/2-\mu_3^2/4$ which gives
$$
P=x^3 y+x^2( 2y^3+\mu_3)+x\left(y^5 + \mu_3 y^2 \right)+ \frac{y^7}7+ \frac{\mu_3 y^4}4 \quad\text{and}
\quad Q=x^2 y+x(y^3+\mu_3)+\frac{y^5}5+\frac{\mu_3 y^2}2.
$$
Note that $P,Q\in K[x,y]$ and $[P,Q]=x^4 y+\mu_3 x^3$. This example is closely related to the example
obtained in~\ref{seccion}, in fact if we apply the procedure of section 1, with $\mu_0=1$,
$\mu_1=0=\mu_2$ and $\mu_3=2$ as in~\ref{seccion}
then we can construct a pair $P,Q\in K[x,y]$ with $\deg(P)=112$, $\deg(Q)=80$ and $[P,Q]=2x^3+x^4 y$.

 The only other solutions were the homogeneous solutions with
$\mu_3=\mu_2=\mu_1=\mu_0=0$.
For $deg(q_1)=5$, after an hour the PC hadn't solved the resulting system.
We also were able to show that in the case $\mu_1=0=\mu_2$ (and $q_1$ with arbitrary degree), any solution of~\eqref{EDPol}
satisfying~\eqref{condiciones} must have $\mu_0=0$.

 Based on this partial results, we state the following conjecture:\\

{\bf CONJECTURE:} The only solutions of~\eqref{EDPol}
are the solutions with $\mu_2=\mu_1=0$.\\

If the conjecture is true, then the only solutions of~\eqref{EDPol} satisfying~\eqref{condiciones}
are the solutions with $\mu_2=\mu_1=\mu_0=0$, which implies $B>16$.

\end{document}